\documentclass[amstex,16pt, amssymb]{article}

\usepackage{mathtext}
\usepackage[cp1251]{inputenc}
\usepackage[T2A]{fontenc}
\usepackage[dvips]{graphicx}
\usepackage{amsmath}
\usepackage{amssymb}
\usepackage{amsxtra}
\usepackage{latexsym}
\usepackage{ifthen}
\usepackage{amsthm}

\numberwithin{equation}{section}

\begin{document}

\title{On boundary value problems for\\  inhomogeneous Beltrami equations}

\author{V. Gutlyanski\u\i{}, O. Nesmelova, \\ V. Ryazanov, E. Yakubov}







\date{}

\maketitle

\begin{abstract}
We give here results on the existence of nonclassical solutions of
the Hilbert boundary value problem in terms of the so-called angular
limits (along nontangent curves to the boundary) for Beltrami
equations with sources in Jordan domains satisfying the
quasihyperbolic boundary condition by Gehring--Martio,
ge\-ne\-ral\-ly speaking, without $(A)-$condition by
Ladyzhenskaya--Ural'tseva and, in particular, without the known
outer cone condition that were standard for boundary value problems
in the PDE theory. Assuming that the coefficients of the problem are
functions of countable bounded variation and the boundary data are
measurable with respect to the logarithmic ca\-pa\-ci\-ty, we prove
the existence of locally H\"older continuous solutions of the
problem.

Moreover, we prove similar results on the Hilbert boundary value
problem with its arbitrary measurable boundary data as well as
coefficients in finitely connected Jordan domains along the
so-called ge\-ne\-ral Bagemihl--Seidel systems of Jordan arcs to
their boundaries. In the same terms, we formulate theorems on the
existence of regular solutions of the known Riemann boundary value
problems, including nonlinear, with arbitrary mea\-su\-rab\-le
coefficients for the nonhomogeneous Beltrami equations. We give also
the representation of obtained solutions through the so--called
generalized analytic functions with sources.

Finally, we formulate similar results on Poincare and Neumann
problems for the Poisson type equations that are main in
hydromechanics (mechanics of in\-com\-pres\-sib\-le fluids) in
ani\-so\-tro\-pic and inhomogeneous media. We give here the
representation of their solutions through the so--called generalized
harmonic functions with sources that describe the corresponding
physical processes in iso\-tro\-pic and homogeneous media.
\end{abstract}

\bigskip



\bigskip



\section{Introduction}


The present paper is a natural continuation of articles \cite{ER2},
\cite{GNR-}--\cite{GRY1}, \cite{Y} and \cite{YR} devoted to the
Riemann, Hilbert, Dirichlet, Poincare and, in par\-ti\-cu\-lar,
Neumann boundary value problems in Jordan domains for
quasiconformal, analytic, harmonic and the so--called
$A-$har\-mo\-nic functions with arbitrary boundary data that are
measurable with respect to lo\-ga\-rith\-mic capacity, as well as
with respect to the natural parameter in domains with rectifiable
boundaries, see \cite{R1}--\cite{R4}. Relevant definitions with
history notes and necessary comments on the previous results can be
found e.g. in \cite{GNRY} and \cite{GRY5}. We recall here only some
of them.


The first part of article \cite{GNRY} was devoted to the proof of
existence of nonclassical solutions of Riemann, Hilbert and
Dirichlet boundary value problems with arbitrary measurable boundary
data with respect to logarithmic capacity for the equations of the
Vekua type
\begin{equation}\label{eqS}
\partial_{\bar z}h(z)\ =\ g(z)
\end{equation}
with real valued functions $g$ in classes $L^{p}(D)$, $p>2$, in the
corresponding domains $D\subset\mathbb C$. We called continuous
solutions $h$ of the equations (\ref{eqS}) with first Sobolev
partial derivatives {\bf generalized analytic functions with the
sources} $\bf g$. Note that these results are valid for complex
va\-lu\-ed sources $g$, see also \cite{GNRY+}, and that equation
(\ref{eqS}) is the complex form of the Poisson equation, which
describes physical processes in homogeneous and isotropic media.


As a consequence, the second part of article \cite{GNRY} contained
the proof of existence of nonclassical solutions to the Poincare
problem on the directional derivatives and, in par\-ti\-cu\-lar, to
the Neumann problem with arbitrary measurable boundary data with
respect to logarithmic ca\-pa\-ci\-ty for the Poisson equations
\begin{equation}\label{eqSSS}
\triangle U(z)\ =\ G(z)
\end{equation}
with real valued functions $G$ of a class $L^{p}(D)$, $p>2$. For
short, we called con\-ti\-nu\-ous solutions to (\ref{eqSSS}) in
$W^{2,p}_{\rm loc}(D)$ {\bf generalized harmonic functions with the
sources } $\bf G$. Note that by the Sobolev embedding theorem, see
Theorem I.10.2 in \cite{So}, such functions belong to the class
$C^1$.

Moreover, in \cite{GRY5} we studied the Hilbert boundary value
problem for the so-called Beltrami equation that is the complex form
of the main equation of the hydromechanics (incompressible fluid
mechanics) in anisotropic and inhomogeneous media, however, without
any sources.


Recall that the {\bf Beltrami equation} is the equation of the form
\begin{equation}\label{1}
f_{\bar{z}}=\mu(z) f_z\
\end{equation}
where $\mu: D\to\mathbb C$ is a mea\-su\-rab\-le function with
$|\mu(z)|<1$ a.e.,  $f_{\bar z}={\bar\partial}f=(f_x+if_y)/2$,
$f_{z}=\partial f=(f_x-if_y)/2$, $z=x+iy$, $f_x$ and $f_y$ are
partial derivatives of the function $f$ in $x$ and $y$,
respectively. Note that continuous functions with the generalized
derivative $f_{\bar z}=0$ are analytic functions, see e.g. Lemma 1
in \cite{ABe}.

Equation~\eqref{1} is said to be {\bf nondegenerate} if
$||\mu||_{\infty}<1$ that we  assume later on. Homeomorphic
solutions $f$ of nondegenerate (\ref{1}) in $W^{1,2}_{\rm loc}$ are
called {\bf quasiconformal mappings} or sometimes {\bf $\mu
-$conformal mappings}. Its continuous solutions in $W^{1,2}_{\rm
loc}$ are called {\bf $\mu -$conformal functions}. Existence
theorems see e.g. in monographs \cite{Alf}, \cite{BGMR} and
\cite{LV}.

\medskip

In the present paper, we study the boundary value problems for the
Beltrami equations with sources. Namely, here we will research the
nonhomogeneous Beltrami equations in the complex plane $\mathbb C$
or in its domains $D$:
\begin{equation}\label{s}
\omega_{\bar{z}}\ =\ \mu(z)\cdot \omega_z\ +\ \sigma (z)\ .
\end{equation}

Following \cite{ABe}, see also monograph \cite{Alf}, let us first
assume that the {\bf source} $\sigma:\mathbb C\to\mathbb C$ belongs
to class $L_p(\mathbb C)$ for some $p>2$ with
\begin{equation}\label{b}
 k\,C_p\ <\ 1\ ,\ \ \ \ k\ :=\ \|\mu\|_{\infty}\ <\ 1\ ,
\end{equation}
where $C_p$ is the norm of the known operator $T:L_p(\mathbb C)\to
L_p(\mathbb C)$ defined through the Cauchy principal limit of the
singular integral
\begin{equation}\label{i}
 (Tg)(\zeta)\ :=\
 \lim\limits_{\varepsilon\to 0}\left\{-\frac{1}{\pi}\int\limits_{|z-\zeta|>\varepsilon}\frac{g(z)}{(z-\zeta)^2}\
 dxdy\right\}\ ,\ \ \ z=x+iy\ .
\end{equation}

As known, $\| Tg\|_2\ =\ \| g\|_2$, i.e. $C_2=1$, and by the Riesz
convexity theorem $C_p\to 1$ as $p\to 2$. Thus, there are such $p$,
whatever the value of $k$ in (\ref{s}).


Let us denote by $B_p$ the Banach space of functions $\omega$,
defined on the whole plane $\mathbb C$, which satisfy a global
H\"older condition of order $1 - 2/p$, which vanish at the origin,
and whose generalized derivatives $\omega_z$ and $\omega_{\bar z}$
exist and belong to $L_p(\mathbb C)$. The norm in $B_p$ is defined
by
\begin{equation}\label{n}
\|\omega\|_{B_p}\ :=\ \sup\limits_{\underset{z_1\ne
z_2}{z_1,z_2\in\mathbb
C,}}\frac{|\omega(z_1)-\omega(z_2)|}{|z_1-z_2|^{1-2/p}}\ +\
\|\omega_z\|_p\ +\ \|\omega_{\bar z}\|_p\ .
\end{equation}

The principal result in \cite{ABe}, Theorem 1, is the following
statement:

\medskip

{\bf Theorem A.}\ {\it Let condition (\ref{b}) hold and $\sigma\in
L_p(\mathbb C)$ for $p>2$. Then the equation (\ref{s}) has a unique
solution $\omega^{\mu,\sigma}\in B_p$. This is the only solution
with $\omega(0) = 0$ and $\omega_z\in L_p(\mathbb C)$.}

\medskip

Its following consequence holds, see Theorem 4 and Lemma 8 in
\cite{ABe}.

\medskip

{\bf Theorem B.}\ {\it Let $\mu :\mathbb C\to\mathbb C$ be in
$L_{\infty}(\mathbb C)$ with compact support and
$k:=\|\mu\|_{\infty}<1$. Then there exists a unique $\mu -$conformal
mapping $f^{\mu}$ in $\mathbb C$ which vanishes at the origin and
satisfies condition $f^{\mu}_z-1\in L_p(\mathbb C)$ for any $p>2$
with (\ref{b}). Moreover, $f^{\mu}(z)=z+\omega^{\mu,\mu}(z)$.}

\section{Factorization of nonhomogeneous equations}

The following simple statement follows by point (vii) of Theo\-rem 5
in \cite{ABe}.


{\bf Remark 1.} Let $D$ be a bounded domain in $\mathbb C$ and let
$\mu : D\to\mathbb C$ be in class $L_{\infty}(D)$ with
$k:=\|\mu\|_{\infty}<1$. If $\omega_1$ and $\omega_2$ are continuous
solutions of (\ref{s}) in $D$ of class $W^{1,2}_{\rm loc}(D)$, then
$\omega_2-\omega_1={\cal A}\circ f^{\mu}$, where $\mu$ was extended
onto $\mathbb C$ by zero outside of $D$ and $\cal A$ is an analytic
function in the domain $D^*:=f^{\mu}(D)$. Note that here we assumed
nothing on $\sigma$.

\medskip

The next representation of solutions of (\ref{s}) is much more
important.

\medskip

{\bf Lemma 1.} {\it Let $D$ be a bounded domain in $\mathbb C$, $\mu
: D\to\mathbb C$ be in class $L_{\infty}(D)$ with
$k:=\|\mu\|_{\infty}<1$ and let $\sigma:D\to\mathbb C$ be in class
$L_p(D)$, $p>2$, with condition (\ref{b}). Then each continuous
solution $\omega$ of equation (\ref{s}) in $D$ of class $W^{1,p}(D)$
has the representation as a composition $h\circ f^{\mu}|_D$, where
$h$ is a generalized analytic function in the domain
$D^*:=f^{\mu}(D)$ with the source $g\in L_{p_*}(D^*)$,
$p_*:=p^2/2(p-1)\in(2,p)$,
\begin{equation}\label{r}
g\ :=\ \left(f^{\mu}_z\cdot\frac{\sigma}{J}\right)\circ
\left(f^{\mu}\right)^{-1}\ ,
\end{equation} where $J$ is the Jacobian of a quasiconformal mapping
$f^{\mu}:\mathbb C\to\mathbb C$ with some extension of $\mu$ onto
$\mathbb C$. Inversely, if $h$ is a ge\-ne\-ra\-li\-zed analytic
function with source (\ref{r}), then $\omega :=h\circ f^{\mu}$ is a
solution of (\ref{s}) of class $C^{\alpha}_{\rm loc}\cap
W^{1,q}_{\rm loc}(D)$, where $\alpha = 1-2/q$ and
$q:=p_*^2/2(p_*-1)\in(2,p_*)$.}


{\bf Remark 2.} Note that if $h$ is a generalized analytic function
with the source $g$ in the domain $D^*$, then $H=h+{\cal A}$ is so
for any analytic function $\cal A$ in $D^*$ but $|{\cal
A}^{\prime}|^p$ can be integrable only locally in $D^*$. The source
in (\ref{r}) is always in class $L_{p_*}(D^*)$,
$p_*:=p^2/2(p-1)\in(2,p)$, in view of Theorem A with $\sigma$
extended onto $\mathbb C$ by zero outside of $D$. Here we may assume
that $\mu$ is extended onto $\mathbb C$ by zero outside of $D$.
However, any other extension of $\mu$ kee\-ping condition (\ref{b})
is suitable here, too.

\medskip

{\bf Proof.} To be short, let us apply here the notation $f$ instead
of $f^{\mu}$. Let us consider the function $h:=\omega\circ f^{-1}$.
First of all, note that by point (iii) of Theorem 5 in \cite{ABe}
$f^*:=f^{-1}|_{D^*}$, $D^*:=f(D)$, is of class $W^{1,p}(D^*)$. Then,
arguing as under the proof of Lemma 10 in \cite{ABe}, we obtain that
$h\in W^{1,p_*}(D^*)$, where $p_*:=p^2/2(p-1)\in(2,p)$. Since
$\omega=h\circ f$, we get also, see e.g. formulas (28) in \cite{ABe}
or formulas I.C(1) in \cite{Alf}, that
$$
\omega_z\ =\ (h_{\zeta}\circ f)\cdot f_z\ +\
(h_{\overline{\zeta}}\circ f)\cdot\overline{f_{\bar z}}\ ,
$$
$$
\omega_{\bar z}\ =\ (h_{\zeta}\circ f)\cdot f_{\bar z}\ +\
(h_{\overline{\zeta}}\circ f)\cdot\overline{f_{z}}\ ,
$$
and, thus,
$$
\sigma(z)\, =\, \omega_{\bar z}\, -\, \mu(z)\omega_z\, =\,
(h_{\overline{\zeta}}\circ f)\overline{f_{z}}(1-|\mu(z)|^2)\, =\,
(h_{\overline{\zeta}}\circ f)J(z)/f_{z}\ ,
$$
where $J(z)=|f_z|^2-|f_{\bar z}|^2=|f_z|^2(1-|\mu(z)|^2)$ is the
Jacobian of $f$, i.e.,
$$
h_{\overline{\zeta}}\ =\ g(\zeta)\ :=
\left(f_z\frac{\sigma}{J}\right)\circ f^{-1}(\zeta)\ .
$$
Similarly, applying Lemma 10 in \cite{ABe} and the Sobolev embedding
theorem, see Theorem I.10.2 in \cite{So}, we come to the inverse
conclusion. \hfill $\Box$

\section{Hilbert problem with angular limits}

In this section, we prove the existence of nonclassical solutions of
the Hilbert boundary value problem with arbitrary boun\-da\-ry data
that are measurable with respect to logarithmic capacity for
nonhomogeneous Beltrami equations. The result is formulated in terms
of the angular limit that is a traditional tool of the geometric
function theory, see e.g. monographs \cite{Du, Ko, L, Po} and
\cite{P}.


Recall that the classic boundary value {\bf problem of Hilbert}, see
\cite{H1}, was formulated as follows: To find an analytic function
$f(z)$ in a domain $D$ bounded by a rectifiable Jordan contour $C$
that satisfies the boundary condition
\begin{equation}\label{2}
\lim\limits_{z\to\zeta , z\in D} \ {\rm{Re}}\,
\{\overline{\lambda(\zeta)}\ f(z)\}\ =\ \varphi(\zeta)
\quad\quad\quad\ \ \ \forall \ \zeta\in C\ ,
\end{equation}
where the {\bf coefficient} $\lambda$ and the {\bf boundary date}
$\varphi$ of the problem are con\-ti\-nu\-ous\-ly differentiable
with respect to the natural parameter $s$ and $\lambda\ne 0$
everywhere on $C$. The latter allows to consider that
$|\lambda|\equiv 1$ on $C$. Note that the quantity
${\rm{Re}}\,\{\overline{\lambda}\, f\}$ in (\ref{2}) means a
projection of $f$ into the direction $\lambda$ interpreted as
vectors in $\mathbb R^2$.

The reader can find a comprehensive treatment of the theory in
excellent books \cite{Be,BW,HKM,TO}. We also recommend to make
familiar with historic surveys in monographs \cite{G,Mus,Ve} on the
topic with an exhaustive bib\-lio\-gra\-phy and take a look at our
recent  papers, see Introduction.

Next, recall that a straight line $L$ is {\bf tangent} to a curve
$\Gamma$ in $\mathbb C$ at a point $z_0\in\Gamma$ if
\begin{equation} \label{eqTANGENT}
\limsup\limits_{z\to z_0, z\in \Gamma}\ \frac{\hbox{dist}\, (z,
L)}{|z-z_0|}\ =\ 0\ .
\end{equation}

Let $D$ be a Jordan domain in $\mathbb C$ with a tangent at a point
$\zeta\in\partial D.$ A path in $D$ terminating at $\zeta$ is called
{\bf nontangential} if its part in a neighborhood of $\zeta$ lies
inside of an angle with the vertex at $\zeta$. The limit along all
nontangential paths at $\zeta$  is called {\bf angular} at the
point.

Following \cite{GRY5}, we say that a Jordan curve $\Gamma$ in
$\mathbb C$ is {\bf almost smooth} if $\Gamma$ has a tangent {\bf
q.e. (quasi everywhere)} with respect to logarithmic capacity, see
e.g. \cite{La} for the term. In particular, $\Gamma$ is almost
smooth if $\Gamma$ has a tangent at all its points except its
countable collection. The nature of such a Jordan curve $\Gamma$ can
be complicated enough because this countable collection can be
everywhere dense in $\Gamma$, see e.g. \cite{DMRV}.


Recall that the {\bf quasihyperbolic distance} between points $z$
and $z_0$ in a domain $D\subset\mathbb C$ is the quantity $$
k_D(z,z_0)\ :=\ \inf\limits_{\gamma} \int\limits_{\gamma}
{ds}/{d(\zeta,\partial D)}\ , $$ where $d(\zeta,\partial D)$ denotes
the Euclidean distance from the point $\zeta\in D$ to $\partial D$
and the infimum is taken over all rectifiable curves $\gamma$
joining the points $z$ and $z_0$ in $D$, see \cite{GP}.

Further, it is said that a domain $D$ satisfies the {\bf
quasihyperbolic boun\-dary con\-di\-tion} if there exist constants $a$
and $b$ and a point $z_0\in D$ such that
\begin{equation} \label{eqHYPERB}
k_D(z,z_0)\ \le\ a\ +\ b\ \ln \frac{d(z_0,\partial D)}{d(z ,\partial
D)}\ \ \ \ \ \ \ \ \ \forall\ z\in D\ .
\end{equation}  The latter notion
was introduced in \cite{GM} but, before it, was first implicitly
applied in \cite{BP}. By the discussion in \cite{GRY5}, every smooth
(or Lipschitz) domain satisfies the qua\-si\-hy\-per\-bo\-lic
boundary condition but such boundaries can be nowhere locally
rectifiable.

Note that it is well--known the so--called $(A)-$condition by
Lady\-zhens\-kaya--Ural'tseva, which is standard in the theory of
boundary value problems for PDE, see e.g. \cite{LU}. Recall that a
domain $D$ in ${\mathbb R}^n$, $n\ge 2$, is called satisfying {\bf
(A)-condition} if
\begin{equation} \label{eqA}
\hbox{mes}\ D\cap B(\zeta,\rho)\ \le\ \Theta_0\,\hbox{mes}\
B(\zeta,\rho)\ \ \ \ \ \ \ \ \forall\ \zeta\in\partial D\ ,\
\rho\le\rho_0
\end{equation}
for some $\Theta_0$ and $\rho_0\in(0,1)$, where $B(\zeta,\rho)$
denotes the ball with the center $\zeta\in{\mathbb R}^n$ and the
radius $\rho$, see 1.1.3 in \cite{LU}.

A domain $D$ in ${\mathbb R}^n$, $n\ge 2$, is said to be satisfying
the {\bf outer cone condition} if there is a cone that makes
possible to be touched by its top to every point of $\partial D$
from the completion of $D$ after its suitable rotations and shifts.
It is clear that the latter condition implies (A)--condition.

Probably one of the simplest examples of an almost smooth domain $D$
with the quasihyperbolic boundary condition and without
(A)--condition is the union of 3 open disks with the radius 1
centered at the points $0$ and $1\pm i$. It is clear that this
domain has zero interior angle at its boundary point $1$.

Given a Jordan domain $D$  in $\mathbb C$, we call $\lambda:\partial
D\to\mathbb C$ a {\bf function of bounded variation}, write
$\lambda\in\mathcal{BV}(\partial D)$, if
\begin{equation}\label{20}
V_{\lambda}(\partial D)\ \colon =\ \sup\
\sum\limits_{j=1}\limits^{k}\
|\lambda(\zeta_{j+1})-\lambda(\zeta_j)| \ <\ \infty \end{equation}
where the supremum is taken over all finite collections of points
$\zeta_j\in\partial D$, $j=1,\ldots , k$, with the cyclic order
meaning that $\zeta_j$ lies between $\zeta_{j+1}$ and $\zeta_{j-1}$
for every $j=1,\ldots , k$. Here we assume that
$\zeta_{k+1}=\zeta_1=\zeta_0$. The quantity $V_{\lambda}(\partial
D)$ is called the {\bf variation of the function} $\lambda$.

Now, we call $\lambda:\partial D\to\mathbb C$ a function of {\bf
countable bounded variation}, write $\lambda\in{\mathcal{
CBV}}(\partial D)$, if there is a countable collection of mutually
disjoint arcs $\gamma_n$ of $\partial D$, $n=1,2,\ldots$ on each of
which the restriction of $\lambda$ is of bounded variation and the
set $\partial D\setminus \cup\gamma_n$ has logarithmic capacity
zero. In particular, the latter holds true if the set $\partial
D\setminus \cup\gamma_n$ is countable. It is clear that such
functions can be singular enough.

\medskip

{\bf Theorem 1.}{\it\, Let $D$ be  a Jordan domain with the
quasihyperbolic boundary condition, $\partial D$ have a tangent
q.e., $\lambda:\partial D\to\mathbb{C},\: |\lambda(\zeta)|\equiv1$,
be in $\mathcal{CBV}(\partial{D})$ and let $\varphi:\partial
D\to\mathbb{R}$ be measurable with respect to logarithmic capacity.
Suppose also that $\mu:D\to\mathbb C$ is of class $L_{\infty}(D)$
with $k:=\|\mu\|_{\infty}<1$, it is H\"older continuous in a
neighborhood of $\partial D$, $\sigma\in L_p(D)$ and condition
(\ref{b}) hold for some $p>2$.

Then equation (\ref{s}) has solutions $\omega: D\to\mathbb C$ of
class $C^{\alpha}_{\rm loc}\cap W^{1,q}_{\rm loc}(D)$, where $\alpha
= 1-2/q$ and $q\in(2,p)$, smooth in the neighborhood of $\partial D$
with the angular limits
\begin{equation}\label{eqLIMH} \lim\limits_{z\to\zeta , z\in D}\ \mathrm
{Re}\ \left\{\, \overline{\lambda(\zeta)}\cdot \omega(z)\, \right\}\
=\ \varphi(\zeta) \quad\quad\quad \mbox{q.e.\ on $\ \partial D$}\ .
\end{equation}
Furthermore, the space of all such solutions $\omega$ of the
equation (\ref{s}) has infinite dimension for each fixed $\lambda$,
$\varphi$, $\mu$ and $\sigma$.}

\medskip

{\bf Remark 3.} By the construction in the proof below, each such
solution has the representation $\omega =h\circ f|_D$, where
$f=f^{\mu}:\mathbb C\to\mathbb C$ is a quasiconformal mapping with
the corresponding extension of $\mu$ onto $\mathbb C$ and $h$ is a
generalized analytic function with the source $g\in L_{p_*}(D_*)$,
$p_*:=p^2/2(p-1)\in(2,p)$, in (\ref{r}), $D_*:=f(D)$, and the
angular limits
\begin{equation}\label{eqLIMH} \lim\limits_{w\to\xi , w\in D_*}\ \mathrm
{Re}\ \left\{\, \overline{\Lambda(\xi)}\cdot h(w)\, \right\}\ =\
\Phi(\xi) \quad\quad\quad \mbox{q.e.\ on $\ \partial D_*$}\ ,
\end{equation} with
$\Lambda\,:=\,\lambda\circ f^{-1}$, $\Phi\,:=\,\varphi\circ f^{-1}$.
Also, $q=p_*^2/2(p_*-1)\in(2,p_*)$.

\bigskip

{\bf Proof.} First of all, let us choose a suitable extension of
$\mu$ onto $\mathbb C$ outside of $D$. By hypotheses of Theorem 1
$\mu$ belongs to a class $C^{\alpha}$, $\alpha\in(0,1)$, for an open
neighborhood $U$ of $\partial D$ inside of $D$. By Lemma 1 in
\cite{GRY1} $\mu$ is extended to a H\"older continuous function
$\mu:U\cup\mathbb C\setminus D\to\mathbb{C}$ of the class
$C^{\alpha}$. Then, for every $k_*\in(k,1)$, there is an open
neighborhood $V$ of $\partial D$ in $\mathbb C$, where
$\|\mu\|_{\infty}\le k_*$ and $\mu$ in $C^{\alpha}(V)$. Let us
choose $k_*\in(k,1)$ so close to $k$ that $k_*C_p<1$ and set
$\mu\equiv 0$ outside of $D\cup V$.

By the Measurable Riemann Mapping Theorem, see e.g. \cite{Alf},
\cite{BGMR} and \cite{LV}, there is a quasiconformal mapping
$f=f^{\mu}:\mathbb C\to{\mathbb{C}}$ a.e. satisfying the Beltrami
equation (\ref{1}) with the given extended complex coefficient $\mu$
in $\mathbb C$. Note that the mapping $f$ has the H\"older
continuous first partial derivatives in $V$ with the same order of
the H\"older continuity as $\mu$, see e.g. \cite{Iw} and also
\cite{IwDis}. Moreover, its Jacobian
\begin{equation}\label{6.13} J(z)\ne 0\ \ \ \ \ \ \ \ \ \ \ \ \
\forall\ z\in V\ , \end{equation} see e.g. Theorem V.7.1 in
\cite{LV}. Hence $f^{-1}$ is also smooth in $V_*:=f(V)$, see e.g.
formulas I.C(3) in \cite{Alf}.

Now, the domain $D_*:=f(D)$ satisfies the boundary quasihyperbolic
condition because $D$ is so, see e.g. Lemma 3.20 in \cite{GM}.
Moreover, $\partial D_*$ has q.e. tangents, furthermore, the points
of $\partial D$ and $\partial D^*$ with tangents correspond each to
other in one-to-one manner because the mappings $f$ and $f^{-1}$ are
smooth there. In addition, the function $\Lambda\,:=\,\lambda\circ
f^{-1}$ belongs to the class $\mathcal{CBV}(\partial{D_*})$ and
$\Phi\,:=\,\varphi\circ f^{-1}$ is measurable with respect to
logarithmic capacity, see e.g. Remark 2.1 in \cite{GRY5}. Next, by
Remark 3 the source $g:D_*\to\mathbb C$ in (\ref{r}) belongs to
class $L_{p_*}(D_*)$, where $p_*=p^2/2(p-1)\in(2,p)$. Thus, by
Theorem 1 in \cite{GNRY} the space of all generalized analytic
functions $h:D_*\to\mathbb C$ with the source $g$ and the angular
limits (\ref{eqLIMH}) q.e. on $\partial D_*$ has infinite dimension.
Finally, by Lemma 1 we obtain the rest of conclusions. \hfill $\Box$

\medskip

In particular case $\lambda\equiv 1$, we obtain the consequence of
Theorem 1 on the Dirichlet problem for the nonhomogeneous Beltrami
equations.

\section{Hilbert problem along special curve systems}

Let $D$ be a domain in $\mathbb C$ whose boundary consists of a
finite collection of mutually disjoint Jordan curves. A family of
mutually disjoint Jordan arcs $J_{\zeta}:[0,1]\to\overline D$,
$\zeta\in\partial D$, with $J_{\zeta}([0,1))\subset D$ and
$J_{\zeta}(1)=\zeta$ that is continuous in the parameter $\zeta$ is
called a {\bf Bagemihl--Seidel system} or, in short, of {\bf class}
${\cal{BS}}$ in $D$.

\medskip

{\bf Lemma 2.} {\it Let $D$ be a bounded domain in $\mathbb C$ whose
boundary consists of a finite number of mutually disjoint Jordan
curves, $\{ \gamma_{\zeta}\}_{\zeta\in\partial D}$ be a family of
Jordan arcs of class ${\cal{BS}}$ and let $\lambda:\partial
D\to\mathbb C$, $|\lambda (\zeta)|\equiv 1$, $\varphi:\partial
D\to\mathbb R$ and $\psi:\partial D\to \mathbb R$ be measurable with
respect to logarithmic capacity.

Suppose also that $\mu:D\to\mathbb C$ is of class $L_{\infty}(D)$
with $k:=\|\mu\|_{\infty}<1$, $\sigma\in L_p(D)$ and condition
(\ref{b}) hold for some $p>2$. Then the equation (\ref{s}) has
solutions $\omega: D\to\mathbb C$ of class $C^{\alpha}_{\rm loc}\cap
W^{1,q}_{\rm loc}$ with $\alpha =1-2/q$ for some $q\in(2,p)$ such
that along the arcs $\gamma_{\zeta}$
\begin{equation}\label{eqARE} \lim\limits_{z\to\zeta}\ \mathrm {Re}\
\{\overline{\lambda(\zeta)}\cdot \omega(z)\}\ =\ \varphi(\zeta)
\quad\quad\quad \mbox{q.e.\ on $\ \partial D$}\ ,
\end{equation}
\begin{equation}\label{eqIM}
\lim\limits_{z\to\zeta}\ \mathrm {Im}\
\{\overline{\lambda(\zeta)}\cdot \omega(z)\}\ =\
\psi(\zeta)\quad\quad\quad \mbox{q.e.\ on $\ \partial D$}\ .
\end{equation}}

\medskip

{\bf Remark 4.} By the construction in the proof below, each such
solution has the representation $\omega =h\circ f|_D$, where
$f=f^{\mu}:\mathbb C\to\mathbb C$ is a quasiconformal mapping with
$\mu$ extended by zero onto $\mathbb C$ outside of $D$ and $h$ is a
generalized analytic function with the source $g\in L_{p_*}(D_*)$,
$p_*:=p^2/2(p-1)\in(2,p)$, in (\ref{r}), $D_*:=f(D)$, and with the
limits along Jordan arcs $\Gamma_{\xi}:=f(\gamma_{f^{-1}(\xi)})$,
$\xi\in\partial D_*$,
\begin{equation}\label{eqLIMHBSRE} \lim\limits_{w\to\xi , w\in D_*}\ \mathrm
{Re}\ \left\{\, \overline{\Lambda(\xi)}\cdot h(w)\, \right\}\ =\
\Phi(\xi) \quad\quad\quad \mbox{q.e.\ on $\ \partial D_*$}\ ,
\end{equation}
\begin{equation}\label{eqLIMHBSIM} \lim\limits_{w\to\xi , w\in D_*}\ \mathrm
{Im}\ \left\{\, \overline{\Lambda(\xi)}\cdot h(w)\, \right\}\ =\
\Psi(\xi) \quad\quad\quad \mbox{q.e.\ on $\ \partial D_*$}\ ,
\end{equation}
$$\mbox{where}\ \ \ \ \Lambda\,:=\,\lambda\circ f^{-1}\ , \Phi\,:=\,\varphi\circ f^{-1}\ ,
\Psi\,:=\,\psi\circ f^{-1}\ .$$ Moreover, more precisely,
$q=p_*^2/2(p_*-1)\in(2,p_*)$.

\bigskip

{\bf Proof.} First of all note that functions $\Lambda$, $\Phi$ and
$\Psi$ in Remark 4 are measurable with respect to logarithmic
capacity on $\partial D_*$ because the quasiconformal mapping $f$ is
H\"older continuous on the compact set $\partial D$, see e.g. Remark
2.1 in \cite{GRY5}. Consequently, by Theorem 2 in \cite{GNRY}, see
also Remark 2 above, there exists a generalized analytic function
$h$ with the source $g$ in (\ref{r}) satisfying conditions
(\ref{eqLIMHBSRE}) and (\ref{eqLIMHBSIM}) q.e. on $\partial D_*$.
Thus, the rest of conclusions follows by Lemma 1, see again Remark
2.1 in \cite{GRY5}, because the inverse quasiconformal mapping
$f^{-1}$ is H\"older continuous on the compact set $\partial D_*$,
too. \hfill $\Box$

\medskip

{\bf Remark 5.} The space of all generalized analytic functions $h$
with the source $g$ in (\ref{r}) satisfying Hilbert boundary
condition (\ref{eqLIMHBSRE}) q.e. on $\partial D_*$ along the Jordan
arcs $\Gamma_{\xi}:=f(\gamma_{f^{-1}(\xi)})$, $\xi\in\partial D_*$,
has infinite dimension for any prescribed $g$, $\Phi$, $\Lambda$ and
$\{ \Gamma_{\xi}\}_{\xi\in\partial D_*}$ of class $\cal BS$ because
the space of all functions $\Psi:\partial D_*\to \mathbb R$ which
are measurable with respect to logarithmic capacity has infinite
dimension.


The latter is valid even for its subspace of continuous functions
$\Psi:\partial D_*\to \mathbb R$. Indeed, every Jordan component of
$\partial D_*$ can be mapped with a homeomorphism onto the unit
circle $\partial\mathbb D$. However, by the Fourier theory, the
space of all continuous functions $\tilde\Psi:\partial\mathbb
D\to\mathbb R$, equivalently, the space of all continuous
$2\pi$-periodic functions $\Psi_*:\mathbb R\to\mathbb R$, just has
infinite dimension.

\medskip

Thus, by Lemma 2, Remarks 4 and 5 we obtain the following statement
on solutions of the equation (\ref{s}), satisfying the Hilbert
boundary condition (\ref{eqARE}) q.e. along Behgemil-Seidel systems
of Jordan arcs.

\medskip

{\bf Theorem 2.}{\it\, Let $D$ be a bounded domain in $\mathbb C$
whose boundary consists of a finite number of mutually disjoint
Jordan curves, $\{ \gamma_{\zeta}\}_{\zeta\in\partial D}$ be a
family of Jordan arcs of class ${\cal{BS}}$ in ${D}$,
$\lambda:\partial D\to\mathbb C$, $|\lambda (\zeta)|\equiv 1$, and
$\varphi:\partial D\to\mathbb R$ be measurable with respect to
logarithmic capacity.

Suppose also that $\mu:D\to\mathbb C$ is of class $L_{\infty}(D)$
with $k:=\|\mu\|_{\infty}<1$, $\sigma\in L_p(D)$ and condition
(\ref{b}) hold for some $p>2$. Then the equation (\ref{s}) has
solutions $\omega: D\to\mathbb C$ of class $C^{\alpha}_{\rm loc}\cap
W^{1,q}_{\rm loc}$ with $\alpha =1-2/q$ for some $q\in(2,p)$ that
satisfy the Hilbert boundary condition (\ref{eqARE}) q.e. in the
sense of the limits along $\gamma_{\zeta}$.

Furthermore, the space of all such solutions $\omega$ has infinite
dimension for any fixed $\mu$, $\sigma$, $\varphi$, $\lambda$ and
$\{ \gamma_{\zeta}\}_{\zeta\in D}$.

Moreover, each such solution has the representation $\omega =h\circ
f|_D$, where $f=f^{\mu}:\mathbb C\to\mathbb C$ is a quasiconformal
mapping with $\mu$ extended by zero onto $\mathbb C$ outside of $D$
and $h$ is a generalized analytic function with the source $g\in
L_{p_*}(D_*)$, $p_*:=p^2/2(p-1)\in(2,p)$, in (\ref{r}), $D_*:=f(D)$,
satisfying the Hilbert boundary condition (\ref{eqLIMHBSRE}) q.e. on
$\partial D_*$ in the sense of the limits along Jordan arcs
$\Gamma_{\xi}:=f(\gamma_{f^{-1}(\xi)})$, $\xi\in\partial D_*$.}

\medskip

In particular case $\lambda\equiv 1$, we obtain the corresponding
consequence on the Dirichlet problem for the nonhomogeneous Beltrami
equations (\ref{s}) along any prescribed Bagemihl--Seidel systems of
Jordan arcs.

\section{Riemann problem and special curve systems}

Recall that the classical setting of the {\bf Riemann problem} in a
smooth Jordan domain $D$ of the complex plane $\mathbb{C}$ is to
find analytic functions $f^+: D\to\mathbb C$ and $f^-:\mathbb
C\setminus \overline{D}\to\mathbb C$ that admit continuous
extensions to $\partial D$ and satisfy the boundary condition
\begin{equation}\label{eqRIEMANN} f^+(\zeta)\ =\
A(\zeta)\cdot f^-(\zeta)\ +\ B(\zeta) \quad\quad\quad \forall\
\zeta\in\partial D
\end{equation}
with its prescribed H\"older continuous coefficients $A:
\partial D\to\mathbb C$ and $B: \partial D\to\mathbb C$.
Recall also that the {\bf Riemann problem with shift} in $D$ is to
find analytic functions $f^+: D\to\mathbb C$ and $f^-:\mathbb
C\setminus \overline{D}\to\mathbb C$ satisfying the condition
\begin{equation}\label{eqSHIFT} f^+(\alpha(\zeta))\ =\ A(\zeta)\cdot
f^-(\zeta)\ +\ B(\zeta) \quad\quad\quad \forall\ \zeta\in\partial D
\end{equation}
where $\alpha :\partial D\to\partial D$ was a one-to-one sense
preserving correspondence having the non-vanishing H\"older
continuous derivative with respect to the natural parameter on
$\partial D$. The function $\alpha$ is called a {\bf shift
function}. The special case $A\equiv 1$ gives the so--called {\bf
jump problem} and then $B\equiv 0$ gives the {\bf problem on gluing}
of analytic functions.

\medskip

Arguing similarly to the proof of Lemma 2, on the base of Lemma 1,
we reduce the proof of the next theorem below to Theorem 3 in
\cite{GNRY} on the Riemann problem for generalized analytic
functions with sources.

\medskip

{\bf Theorem 3.}{\it\, Let $D$ be a domain in ${\mathbb C}$ whose
boundary consists of a finite number of mutually disjoint Jordan
curves, functions $A: \partial D\to\mathbb C$ and $B:
\partial D\to\mathbb C$ be measurable with respect to
logarithmic capacity, $\{\gamma^+_{\zeta}\}_{\zeta\in\partial D}$
and $\{\gamma^-_{\zeta}\}_{\zeta\in\partial D}$ be families of
Jordan arcs of class ${\cal{BS}}$ in ${D}$ and $\mathbb
C\setminus\overline{ D}$, correspondingly. Suppose also that
$\mu:\mathbb C\to\mathbb C$ is of class $L_{\infty}(\mathbb C)$ with
$k:=\|\mu\|_{\infty}<1$, $\sigma\in L_p(\mathbb C)$ has compact
support and condition (\ref{b}) hold for some $p>2$.

Then the nonhomogeneous Beltrami equation (\ref{s}) has solutions
$\omega^+: D\to\mathbb C$ and $\omega^-:{\mathbb
C}\setminus\overline{D}\to\mathbb C$ of class $C^{\alpha}_{\rm
loc}\cap W^{1,q}_{\rm loc}$ with $\alpha =1-2/q$ for some
$q\in(2,p)$, satisfying the Riemann boundary condition
(\ref{eqRIEMANN}) q.e. on $\partial D$, where $\omega^+(\zeta)$ and
$\omega^-(\zeta)$ are limits of $\omega^+(z)$ and $\omega^-(z)$ az
$z\to\zeta$ along $\gamma^+_{\zeta}$ and $ \gamma^-_{\zeta}$,
correspondingly.

Furthermore, the space of all such couples of solutions
$(\omega^+,\omega^-)$ of (\ref{s}) has infinite dimension for any
fixed $\mu$, $\sigma$, couples $(A, B)$ and collections
$\gamma^+_{\zeta}$ and $ \gamma^-_{\zeta}$, $\zeta\in\partial D$.}

\medskip

Moreover, each such couples of solutions $(\omega^+,\omega^-)$ of
(\ref{s}) has the representation in the form of the composition of
the corresponding couples $(h^+,h^-)$ of ge\-ne\-ra\-li\-zed
ana\-ly\-tic functions with the source $g$ in (\ref{r}) and
$f^{\mu}:\mathbb C\to\mathbb C$.

\medskip

Theorem 3 is a special case of the following lemma, whose proof is
reduced to Lemma 1 in \cite{GNRY} on the Riemann boundary value
problem with shifts for ge\-ne\-ra\-li\-zed analytic functions with
source $g$ given by (\ref{r}) on the base of Lemma 1 above, that may
have of independent interest.

\medskip

{\bf Lemma 3.} {\it\, Under the hypotheses of Theorem 3, let in
addition that $\alpha : \partial D\to\partial D$ be a homeomorphism
keeping components of $\partial D$ such that $\alpha$ and
$\alpha^{-1}$ have $(N)-$property with respect to logarithmic
capacity.

Then the nonhomogeneous Beltrami equation (\ref{s}) has solutions
$\omega^+: D\to\mathbb C$ and $\omega^-:{\mathbb
C}\setminus\overline{D}\to\mathbb C$ of class $C^{\alpha}_{\rm
loc}\cap W^{1,q}_{\rm loc}$ with $\alpha =1-2/q$ for some
$q\in(2,p)$, satisfying the Riemann boundary condition with shift
(\ref{eqSHIFT}) q.e. on $\partial D$, where $\omega^+(\zeta)$ and
$\omega^-(\zeta)$ are limits of $\omega^+(z)$ and $\omega^-(z)$ az
$z\to\zeta$ along $\gamma^+_{\zeta}$ and $ \gamma^-_{\zeta}$,
correspondingly.

Furthermore, the space of all such couples of solutions
$(\omega^+,\omega^-)$ has infinite dimension for any fixed $\mu$,
$\sigma$, couples $(A, B)$ and collections $\gamma^+_{\zeta}$ and $
\gamma^-_{\zeta}$, $\zeta\in\partial D$.}

\medskip

Again, each such couples of solutions $(\omega^+,\omega^-)$ of
(\ref{s}) has the representation in the form of the composition of
the corresponding couples $(h^+,h^-)$ of ge\-ne\-ra\-li\-zed
ana\-ly\-tic functions with the source $g$ in (\ref{r}) and
$f^{\mu}:\mathbb C\to\mathbb C$.

\section{Nonlinear Riemann boundary value problems}

We are able by our scheme above to formulate also a series of
results on nonlinear Riemann boundary value problems in terms of
Bagemihl--Seidel systems for Beltrami equations with sources and
representations.

For instance, special nonlinear boundary value problems of the form
\begin{equation}\label{eqNONLINEAR} \omega^+(\zeta)\ =\ \varphi(\,\zeta,\,
 \omega^-(\zeta)\, ) \quad\quad\quad \mbox{q.e. on}\quad \zeta\in\partial D
\end{equation}
are solved if $\varphi : \partial D\times\mathbb C\to\mathbb C$
satisfies the {\bf Cara\-theo\-dory conditions} with respect to
logarithmic capacity, i.e., if $\varphi(\zeta, w)$ is continuous in
the variable $w\in\mathbb C$ for q.e. $\zeta\in\partial D$ and it is
measurable with respect to logarithmic capacity in the variable
$\zeta\in\partial D$ for all $w\in\mathbb C$. Later on, we sometimes
say in short "C-measurable" instead of the expression "measurable
with respect to logarithmic capacity".

\medskip

{\bf Theorem 4.}{\it\, Let $D$ be a domain in ${\mathbb C}$ whose
boundary consists of a finite number of mutually disjoint Jordan
curves, functions $A: \partial D\to\mathbb C$ and $B:
\partial D\to\mathbb C$ be measurable and $\varphi : \partial D\times\mathbb C\to\mathbb C$
satisfy the Cara\-theo\-dory conditions with respect to logarithmic
capacity, $\{\gamma^+_{\zeta}\}_{\zeta\in\partial D}$ and
$\{\gamma^-_{\zeta}\}_{\zeta\in\partial D}$ be families of Jordan
arcs of class ${\cal{BS}}$ in ${D}$ and $\mathbb
C\setminus\overline{ D}$, correspondingly.

Suppose also that $\mu:\mathbb C\to\mathbb C$ is of class
$L_{\infty}(\mathbb C)$ with $k:=\|\mu\|_{\infty}<1$, $\sigma\in
L_p(\mathbb C)$ has compact support and condition (\ref{b}) hold for
some $p>2$. Then the nonhomogeneous Beltrami equation (\ref{s}) has
solutions $\omega^+: D\to\mathbb C$ and $\omega^-:{\mathbb
C}\setminus\overline{D}\to\mathbb C$ of class $C^{\alpha}_{\rm
loc}\cap W^{1,q}_{\rm loc}$ with $\alpha =1-2/q$ with some
$q\in(2,p)$, satisfying the nonlinear Riemann boundary condition
(\ref{eqNONLINEAR}) q.e. on $\partial D$, where $\omega^+(\zeta)$
and $\omega^-(\zeta)$ are limits of $\omega^+(z)$ and $\omega^-(z)$
az $z\to\zeta$ along $\gamma^+_{\zeta}$ and $ \gamma^-_{\zeta}$,
correspondingly.

Furthermore, the space of all such couples of solutions
$(\omega^+,\omega^-)$ has infinite dimension for any fixed $\mu$,
$\sigma$, couples $(A, B)$ and collections $\gamma^+_{\zeta}$ and $
\gamma^-_{\zeta}$, $\zeta\in\partial D$.}

\medskip

{\bf Proof.} The spaces of solutions of such problems always have
infinite dimension. Indeed, by the Egorov theorem, see e.g. Theorem
2.3.7 in \cite{Fe}, see also Section 17.1 in \cite{KZPS}, the
function $\varphi(\zeta,\psi(\zeta))$ is $C-$measurable in
$\zeta\in\partial D$ for every $C-$measurable function
$\psi:\partial D\to\mathbb C$ if the function $\varphi$ satisfies
the {Caratheodory conditions}, and the space of all $C-$measurable
functions $\psi:\partial D\to\mathbb C$ has the infinite dimension,
see e.g. arguments in Remark 5.

Finally, applying Theorem 3 firstly to each component of $\mathbb
C\setminus\overline D$ with $A\equiv 0$ and $B=\psi$ on the
corresponding component of $\partial D$, we see that $\omega^-$ can
be just $\psi$, and secondly to the domain $D$ with $A\equiv 0$ and
$B(\zeta):=\varphi(\zeta,\psi(\zeta))$, $\zeta\in\partial D$, we
come to the desired conclusion. \hfill $\Box$

\medskip

Similarly we also obtain the following consequence of Lemma 3 on the
nonlinear Riemann boundary value problems with shifts.

\medskip

{\bf Lemma 4.} {\it\, Under the hypotheses of Theorem 4, let in
addition that $\alpha : \partial D\to\partial D$ be a homeomorphism
keeping components of $\partial D$ such that $\alpha$ and
$\alpha^{-1}$ have $(N)-$property with respect to logarithmic
capacity.

Then the nonhomogeneous Beltrami equation (\ref{s}) has solutions
$\omega^+: D\to\mathbb C$ and $\omega^-:{\mathbb
C}\setminus\overline{D}\to\mathbb C$ of class $C^{\alpha}_{\rm
loc}\cap W^{1,q}_{\rm loc}$ with $\alpha =1-2/q$ for some
$q\in(2,p)$, satisfying the nonlinear Riemann boundary condition
with shift
\begin{equation}\label{eqNONLINEARSHIFT} \omega^+(\alpha(\zeta))\ =\ \varphi(\,\zeta,\,
 \omega^-(\zeta)\, ) \quad\quad\quad \mbox{q.e. on}\quad \zeta\in\partial
 D\ ,
\end{equation}
where $\omega^+(\zeta)$ and $\omega^-(\zeta)$ are limits of
$\omega^+(z)$ and $\omega^-(z)$ az $z\to\zeta$ along
$\gamma^+_{\zeta}$ and $ \gamma^-_{\zeta}$, correspondingly.

Furthermore, the space of all such couples of solutions
$(\omega^+,\omega^-)$ has infinite dimension for any fixed $\mu$,
$\sigma$, couples $(A, B)$ and collections $\gamma^+_{\zeta}$ and $
\gamma^-_{\zeta}$, $\zeta\in\partial D$.}

\medskip

Representations of solutions in Theorem 4 and Lemma 4 as
compositions of generalized analytic functions $h^+$ and $h^-$ with
the source $g$ in (\ref{r}) and quasiconformal mappings
$f^{\mu}:\mathbb C\to\mathbb C$ remain also valid.

\section{Mixed nonlinear boundary-value problems}

In order to demonstrate the potentiality of our approach, we give
here also a couple of results on mixed nonlinear boundary value
problems in terms of Bagemihl--Seidel systems for Beltrami equations
with sources.

\medskip

Namely, arguing similarly to the last section and to the proof of
Theorem 1, we are able to reduce the proof of the following theorem
to Theorem 4 in \cite{GNRY} on the base of Lemma 1 and Remark 2.

\medskip

{\bf Theorem 5.}{\it\, Let $D$ be a domain in ${\mathbb C}$ whose
boundary consists of a finite number of mutually disjoint Jordan
curves, $\varphi :
\partial D\times\mathbb C\to\mathbb C$ satisfy the Carath\'{e}odory
conditions and $\nu:\partial D\to\mathbb C$, $|\nu (\zeta)|\equiv
1$, be measurable with respect to the logarithmic capacity,
$\{\gamma^+_{\zeta}\}_{\zeta\in\partial D}$ and
$\{\gamma^-_{\zeta}\}_{\zeta\in\partial D}$ be families of Jordan
arcs of class ${\cal{BS}}$ in ${D}$ and $\mathbb
C\setminus\overline{ D}$, correspondingly.

Suppose also that $\mu:\mathbb C\to\mathbb C$ is of class
$L_{\infty}(\mathbb C)$ with $k:=\|\mu\|_{\infty}<1$, $\sigma\in
L_p(\mathbb C)$ has compact support and condition (\ref{b}) hold for
some $p>2$, they are both H\"older continuous in an open
neighborhood $V$ of $\partial D$.

Then nonhomogeneous Beltrami equation (\ref{s}) has solutions
$\omega^+: D\to\mathbb C$ and $\omega^-:{\mathbb
C}\setminus\overline{D}\to\mathbb C$ of class $C^{\alpha}_{\rm
loc}\cap W^{1,q}_{\rm loc}$ with $\alpha =1-2/q$ for some
$q\in(2,p)$, that both are smooth in $V$ and that satisfy the mixed
nonlinear boundary condition
\begin{equation}\label{eqMIXED} \omega^+(\zeta)\ =\ \varphi\left(\,\zeta,\,
\left[\frac{\partial \omega}{\partial\nu}\right]^- (\zeta)\, \right)
\quad\quad\quad \mbox{q.e.\ on $\ \partial D$}\ ,
\end{equation}
where $\omega^+(\zeta)$ and $\left[\frac{\partial
\omega}{\partial\nu}\right]^- (\zeta)$ are limits of the functions
$\omega^+(z)$ and $\frac{\partial \omega^-}{\partial\, \nu\,}\ (z)$
as $z\to\zeta$ along $\gamma^+_{\zeta}$ and $ \gamma^-_{\zeta}$,
correspondingly.


Furthermore, the space of all such couples $(\omega^+,\omega^-)$ has
infinite dimension for any such prescribed functions $\sigma$,
$\varphi$, $\nu$ and collections $\gamma^+_{\zeta}$ and $
\gamma^-_{\zeta}$, $\zeta\in\partial D$.}

\medskip

Indeed, by the Measurable Riemann Mapping Theorem, see e.g.
\cite{Alf}, \cite{BGMR} and \cite{LV}, there is a quasiconformal
mapping $f=f^{\mu}:\mathbb C\to{\mathbb{C}}$ a.e. satisfying the
Beltrami equation (\ref{1}) with the given complex coefficient $\mu$
in $\mathbb C$. Note also that the mapping $f$ has locally H\"older
continuous first partial derivatives in $V$ with the same order of
the H\"older continuity as $\mu$, see e.g. \cite{Iw} and also
\cite{IwDis}. Moreover, its Jacobian
\begin{equation}\label{6.13} J(z)\ne 0\ \ \ \ \ \ \ \ \ \ \ \ \
\forall\ z\in V\ , \end{equation} see e.g. Theorem V.7.1 in
\cite{LV}. Hence $f^{-1}$ is also smooth in $V_*:=f(V)$, see e.g.
formulas I.C(3) in \cite{Alf}. Moreover, by Lemma 10 in \cite{ABe}
for $\omega:=\omega^-$ and $h=h^-$ we have the connection between
derivatives of $\omega$ and $h$ in the corresponding directions
$$
\frac{\partial\omega}{\partial\nu}\ =\ \omega_z\cdot \nu\ +\
\omega_{\bar z}\cdot\overline\nu\ =\ \nu\cdot(h_w\circ f\cdot f_z\
+\ h_{\bar w}\circ f\cdot\overline{f_{\bar z}})\ +\
\overline{\nu}\cdot(h_w\circ f\cdot f_z\ +\ h_{\bar w}\circ
f\cdot\overline{f_{z}})
$$
$$
=\ h_w\circ f\cdot (\nu\cdot f_{z}\ +\ \overline{\nu}\cdot f_{\bar
z})\
 +\ h_{\overline w}\circ f\cdot (\nu\cdot \overline{f_{\bar z}}\ +\ \overline{\nu}\cdot
 \overline{f_{z}})\ =\ h_w\circ f\cdot \frac{\partial f}{\partial
 \nu}\ +\ h_{\bar w}\circ f\cdot\overline{\frac{\partial f}{\partial\nu}}
$$
$$
=\ \left(\ {\cal N}_*\cdot h_w\ +\ \overline{\cal N}_*\cdot h_{\bar
w}\ \right)\circ f\ =\ \frac{\partial h}{\partial{\cal N}_*}\circ f\
,\ \ \ \ \ \ {\cal N}_*\ :=\ \frac{\partial f}{\partial \nu}\circ
f^{-1}\ .
$$

\medskip

Theorem 5 is a special case of the following lemma on the
corresponding mixed boundary value problem with shift.

\medskip

{\bf Lemma 5.}{\it\, Under the hypotheses of Theorem 4, let in
addition $\beta : \partial D\to\partial D$ be a homeomorphism
keeping components of $\partial D$ such that $\beta$ and
$\beta^{-1}$ have $(N)-$property with respect to logarithmic
capacity.

Then nonhomogeneous Beltrami equation (\ref{s}) has solutions
$\omega^+: D\to\mathbb C$ and $\omega^-:{\mathbb
C}\setminus\overline{D}\to\mathbb C$ of class $C^{\alpha}_{\rm
loc}\cap W^{1,q}_{\rm loc}$ with $\alpha =1-2/q$ for some
$q\in(2,p)$, that both are smooth about $\partial D$ and that
satisfy the following mixed nonlinear boundary condition
\begin{equation}\label{eqMIXED_SHIFT} \omega^+(\beta(\zeta))\ =\
\varphi\left(\,\zeta,\, \left[\frac{\partial
\omega}{\partial\nu}\right]^- (\zeta)\, \right) \quad\quad\quad
\mbox{q.e.\ on $\ \partial D$}\ ,
\end{equation}
where $\omega^+(\zeta)$ and $\left[\frac{\partial
\omega}{\partial\nu}\right]^- (\zeta)$ are limits of the functions
$\omega^+(z)$ and $\frac{\partial \omega^-}{\partial\, \nu\,}\ (z)$
as $z\to\zeta$ along $\gamma^+_{\zeta}$ and $ \gamma^-_{\zeta}$,
correspondingly.

Furthermore, the space of all such couples $(\omega^+,\omega^-)$ has
infinite dimension for any such prescribed functions $\sigma$,
$\varphi$, $\nu$ and collections $\gamma^+_{\zeta}$ and $
\gamma^-_{\zeta}$, $\zeta\in\partial D$. }

\medskip

{\bf Remark 6.} On the base of Remark 6 in \cite{GNRY}, under the
hypotheses of Theorem 5,
 the similar statement can be derived for
the boundary gluing conditions of the form
\begin{equation}\label{eqTWICE} \left[\frac{\partial
\omega^+}{\partial\nu_*}\right] (\beta(\zeta))\ =\
\varphi\left(\,\zeta,\, \left[\frac{\partial
\omega^-}{\partial\nu}\right] (\zeta)\, \right)\ \ \ \ \ \mbox{q.e.
on $\partial D$}\ .
\end{equation}

\medskip

Again, representations of solutions in Theorem 5, Lemma 5 and Remark
6 as compositions of generalized analytic functions $h^+$ and $h^-$
with the source $g$ in (\ref{r}) and quasiconformal mappings
$f^{\mu}:\mathbb C\to\mathbb C$ remain valid.

\section{Poincar\'{e} and Neumann problems in terms of angular limits}

It is well--known, see Theorem 16.1.6 in \cite{AIM}, that
nonhomogeneous Beltrami equations in a domain $D$ of the complex
plane $\mathbb C$ are closely connected with the divergence type
equations of the form
\begin{equation}\label{6.8} {\rm div}\,
A(z)\,\nabla\,u(z)\ =\ g(z)\ , \end{equation} where $A(z)$ is the
matrix function:
\begin{equation}\label{6.9}
A=\left(\begin{array}{ccc} {|1-\mu|^2\over 1-|\mu|^2}  & {-2{\rm Im}\,\mu\over 1-|\mu|^2} \\
                            {-2{\rm Im}\,\mu\over 1-|\mu|^2}          & {|1+\mu|^2\over 1-|\mu|^2}  \end{array}\right).
\end{equation}
As we see, the matrix function $A(z)$ in (\ref{6.9}) is symmetric
and its entries $a_{ij}=a_{ij}(z)$ are dominated by the quantity
$$
K_{\mu}(z)\ :=\ \frac{1\ +\ |\mu(z)|}{1\ -\ |\mu(z)|}\ ,
$$
and, thus, they are bounded if the Beltrami equation is not
de\-ge\-ne\-ra\-te.


Vice verse, uniformly elliptic equations (\ref{6.8}) with symmetric
$A(z)$ whose entries are measurable and ${\rm det}\,A(z)\equiv 1$
just correspond to nondegenerate Beltrami equations with coefficient
\begin{equation}\label{6.10}
\mu\ =\ \frac{1}{\mathrm{det}\, (I+A)}\ (a_{22}-a_{11}\ -\
2ia_{21})\ =\ \frac{a_{22}-a_{11}\ -\ 2ia_{21}}{1\ +\ \mathrm{Tr}\,
A\ +\ \mathrm{det}\, A}\ .\end{equation} $M^{2\times 2}(D)$ denotes
the collection of all such matrix functions $A(z)$ in $D$.

Note that (\ref{6.8}) are the main equation of hydromechanics
(mechanics of incompressible fluids) in anisotropic and
inhomogeneous media.

Given such a matrix function $A$ and a quasiconformal mapping
$f^{\mu}:D\to\mathbb C$, we have already seen in Lemma 1 of
\cite{GNR-}, by direct computation, that if a function $T$ and the
entries of $A$ are sufficiently smooth, then
\begin{equation}\label{6.10k}
\mbox{div}\,[A(z)\cdot\nabla\,(T(f^{\mu}(z)))]\ =\ J(z)\cdot
\triangle\,T(f^{\mu}(z))\ . \end{equation} In the case $T\in
W^{1,2}_{\rm loc}$, we understand the identity (\ref{6.10k}) in the
distributional sense, see Proposition 3.1 in \cite{GNR}, i.e., for
all $\psi\in C^1_0(D)$,
\begin{equation}\label{6.10kd}
\int\limits_D\langle A\nabla(T\circ f^{\mu}),\nabla\psi\rangle\
dm_z\ =\ \int\limits_DJ\cdot\langle M^{-1}((\nabla T)\circ
f^{\mu}),\nabla\psi\rangle\ dm_z\ ,
\end{equation}
where $M$ is the Jacobian matrix of the mapping $f^{\mu}$ and $J$ is
its Jacobian.

\medskip

{\bf Theorem 6.}{\it\, Let $D$ be  a Jordan domain with the
quasihyperbolic boundary condition, $\partial D$ have a tangent
q.e., $\nu:\partial D\to\mathbb{C},\: |\nu(\zeta)|\equiv 1$, be in
$\mathcal{CBV}(\partial{D})$ and $\varphi:\partial D\to\mathbb{R}$
be measurable with respect to logarithmic capacity. Suppose that
$A\in M^{2\times 2}(D)$ has H\"older continuous entries and $g\in
L_p(D)$ for some $p>2 $.

Then there exist smooth weak solutions $u:D\to\mathbb R$ of equation
(\ref{6.8}) that have the angular limits of its derivatives in
directions $\nu$
\begin{equation}\label{eqLIMIT} \lim\limits_{z\to\zeta , z\in D}\ \frac{\partial u}{\partial \nu}\ (z)\ =\
\varphi(\zeta) \quad\quad\quad \mbox{q.e.\ on $\ \partial D$}\
.\end{equation} Furthermore, the space of all such solutions has
infinite dimension for each fixed collection of functions $A$, $g$,
$\nu$ and $\varphi$ of given classes.}

\medskip

Here $u$ is called a weak solution of equation (\ref{6.8}) if
\begin{equation}\label{W}
\int\limits_{D}\{\langle A(z)\nabla u(z),\nabla \psi\rangle\ +\
g(z)\cdot\psi(z)\}\ dm(z)\ =\ 0\ \ \ \ \forall\ \psi\in C^1_0(D)\ .
\end{equation}


{\bf Remark 7.} By the construction in the proof below, each such
solution $u$ has the representation $u =U\circ f|_D$, where
$f=f^{\mu}:\mathbb C\to\mathbb C$ is a quasiconformal mapping with a
suitable extension of $\mu$ in (\ref{6.10}) and $U$ is a generalized
harmonic function with the source ($J$ is Jacobian of $f$)
\begin{equation}\label{SOURCE}
G:=(g/J)\circ f^{-1} \end{equation} of class $L_{p}(D_*)$ in the
domain $D_*:=f(D)$ that has the angular limits
\begin{equation}\label{eqLIMHharm} \lim\limits_{w\to\xi , w\in D_*}\ \frac{\partial U}{\partial {\cal N}}\ (w)\ =\
\Phi(\xi) \quad\quad\quad \mbox{q.e.\ on $\ \partial D_*$}\ ,
\end{equation} with
\begin{equation}\label{N}
{\cal N}(\xi)\ :=\ \left\{\frac{\partial f}{\partial
\nu}\cdot\left|\frac{\partial f}{\partial
\nu}\right|^{-1}\right\}\circ f^{-1}(\xi)\ ,\ \ \ \xi\in\partial
D_*\ ,
\end{equation}
and
\begin{equation}\label{F}
\Phi(\xi)\ :=\ \left\{\varphi\cdot\left|\frac{\partial f}{\partial
\nu}\right|^{-1}\right\}\circ f^{-1}(\xi)\ ,\ \ \ \xi\in\partial
D_*\ .
\end{equation}


{\bf Proof.} By the hypotheses of the theorem $\mu$ given by
(\ref{6.10}) belongs to a class $C^{\alpha}(D)$, $\alpha\in(0,1)$,
and by Lemma 1 in \cite{GRY1} $\mu$ is extended to a H\"older
continuous function $\mu:\mathbb C\to\mathbb{C}$ of the class
$C^{\alpha}$. Then, for every $k_*\in(k,1)$, there is an open
neighborhood $V$ of $\overline D$, where $\|\mu\|_{\infty}\le k_*$
and $\mu$ is of class $C^{\alpha}(V)$. We set $\mu\equiv 0$ in
$\mathbb C\setminus V$.

By the Measurable Riemann Mapping Theorem, see e.g. \cite{Alf},
\cite{BGMR} and \cite{LV}, there is a quasiconformal mapping
$f=f^{\mu}:\mathbb C\to{\mathbb{C}}$ a.e. satisfying the Beltrami
equation (\ref{1}) with the given extended complex coefficient $\mu$
in $\mathbb C$. Note that the mapping $f$ has the H\"older
continuous first partial derivatives in $V$ with the same order of
the H\"older continuity as $\mu$, see e.g. \cite{Iw} and also
\cite{IwDis}. Moreover, its Jacobian
\begin{equation}\label{6.13} J(z)\ne 0\ \ \ \ \ \ \ \ \ \ \ \ \
\forall\ z\in V\ , \end{equation} see e.g. Theorem V.7.1 in
\cite{LV}. Hence $f^{-1}$ is also smooth in $V_*:=f(V)$, see e.g.
formulas I.C(3) in \cite{Alf}.

Now, the domain $D_*:=f(D)$ satisfies the boundary quasihyperbolic
condition because $D$ is so, see e.g. Lemma 3.20 in \cite{GM}.
Moreover, $\partial D_*$ has q.e. tangents, furthermore, the points
of $\partial D$ and $\partial D^*$ with tangents correspond each to
other in a one-to-one manner because the mappings $f$ and $f^{-1}$
are smooth. In addition, the function $\cal N$ in (\ref{N}) belongs
to the class $\mathcal{CBV}(\partial{D_*})$ and $\Phi$ in (\ref{F})
is measurable with respect to logarithmic capacity, see e.g. Remark
2.1 in \cite{GRY5}. Next, the source $G:D_*\to\mathbb R$ in
(\ref{SOURCE}) belongs to class $L_{p}(D_*)$, see e.g. point (vi) of
Theorem 5 in \cite{ABe}. Thus, by Theorem 5 in \cite{GNRY} the space
of all generalized analytic functions $U:D_*\to\mathbb R$ with the
source $G$ and the angular limits (\ref{eqLIMHharm}) q.e. on
$\partial D_*$ has infinite dimension. Finally, by Proposition 3.1
in \cite{GNR} the functions $u:=U\circ f$ give the desired solutions
of equation (\ref{6.8}) because by Lemma 10 in \cite{ABe}
$$
\frac{\partial u}{\partial\nu}\ =\ u_z\cdot \nu\ +\ u_{\bar
z}\cdot\overline\nu\ =\ \nu\cdot(U_w\circ f\cdot f_z\ +\ U_{\bar
w}\circ f\cdot\overline{f_{\bar z}})\ +\
\overline{\nu}\cdot(U_w\circ f\cdot f_z\ +\ U_{\bar w}\circ
f\cdot\overline{f_{z}})
$$
$$
=\ U_w\circ f\cdot (\nu\cdot f_{z}\ +\ \overline{\nu}\cdot f_{\bar
z})\
 +\ U_{\overline w}\circ f\cdot (\nu\cdot \overline{f_{\bar z}}\ +\ \overline{\nu}\cdot
 \overline{f_{z}})\ =\ U_w\circ f\cdot \frac{\partial f}{\partial
 \nu}\ +\ U_{\bar w}\circ f\cdot\overline{\frac{\partial f}{\partial\nu}}
$$
$$
=\ \left(\ {\cal N}\cdot U_w\ +\ \overline{\cal N}\cdot U_{\bar w}\
\right)\circ f\cdot \left| \frac{\partial f}{\partial \nu}\right|\
=\ \frac{\partial U}{\partial{\cal N}}\circ f\cdot \left|
\frac{\partial f}{\partial \nu}\right|\ ,
$$
where the direction $\cal N$ is given by (\ref{N}).\hfill $\Box$

\bigskip

{\bf Remark 8.} We are able to say more in the case of $\mathrm
{Re}\ n(\zeta)\overline{\nu(\zeta)}>0$, where $n(\zeta)$ is  the
inner normal to $\partial D$ at the point $\zeta$. Indeed, the
latter magnitude is a scalar product of $n=n(\zeta)$ and $\nu
=\nu(\zeta)$ interpreted as vectors in $\mathbb R^2$ and it has the
geometric sense of projection of the vector $\nu$ into $n$. In view
of (\ref{eqLIMIT}), since the limit $\varphi(\zeta)$ is finite,
there is a finite limit $u(\zeta)$ of $u(z)$ as $z\to\zeta$ in $D$
along the straight line passing through the point $\zeta$ and being
parallel to the vector $\nu$ because along this line
\begin{equation}\label{eqDIFFERENCE} u(z)\ =\ u(z_0)\ -\ \int\limits_{0}\limits^{1}\
\frac{\partial u}{\partial \nu}\ (z_0+\tau (z-z_0))\ d\tau\
.\end{equation} Thus, at each point with condition (\ref{eqLIMIT}),
there is the directional derivative
\begin{equation}\label{eqPOSITIVE}
\frac{\partial u}{\partial \nu}\ (\zeta)\ :=\ \lim_{t\to 0}\
\frac{u(\zeta+t\cdot\nu)-u(\zeta)}{t}\ =\ \varphi(\zeta)\ .
\end{equation}

\bigskip

In particular, in the case of the Neumann problem, $\mathrm {Re}\
n(\zeta)\overline{\nu(\zeta)}\equiv 1>0$, where $n=n(\zeta)$ denotes
the unit interior normal to $\partial D$ at the point $\zeta$, and
we have by Theorem 6 and Remark 8 the following significant result.

\bigskip

{\bf Corollary 1.} {\it Let $D$ be a Jordan domain in $\Bbb C$ with
the quasihyperbolic boundary condition, the unit inner normal
$n(\zeta)$, $\zeta\in\partial D$, belong to the class
$\mathcal{CBV}(\partial D)$ and $\varphi:\partial D\to\mathbb{R}$ be
measurable with respect to logarithmic capacity. Suppose that $A\in
M^{2\times 2}(D)$ has H\"older continuous entries and $g\in L_p(D)$
for some $p>2 $.

\medskip

Then there exist smooth weak solutions $u:D\to\mathbb R$ of equation
(\ref{6.8}) such that q.e. on $\partial D$ there exist:

\bigskip

1) the finite limit along the normal $n(\zeta)$
$$
u(\zeta)\ :=\ \lim\limits_{z\to\zeta}\ u(z)\ ,$$

2) the normal derivative
$$
\frac{\partial u}{\partial n}\, (\zeta)\ :=\ \lim_{t\to 0}\
\frac{u(\zeta+t\cdot n(\zeta))-u(\zeta)}{t}\ =\ \varphi(\zeta)\ ,
$$

3) the angular limit
$$ \lim_{z\to\zeta}\ \frac{\partial u}{\partial n}\, (z)\ =\
\frac{\partial u}{\partial n}\, (\zeta)\ .$$

Furthermore, the space of such solutions $u$ of has infinite
dimension.}

\bigskip

{\bf Remark 9.} Moreover, such solutions $u$ have the representation
$u =U\circ f|_D$, where $f=f^{\mu}:\mathbb C\to\mathbb C$ is a
quasiconformal mapping with a suitable extension of $\mu$ in
(\ref{6.10}) onto $\mathbb C$ outside of $D$ described in the proof
of Theorem 6 and $U$ is a generalized harmonic function with the
source $G$ in (\ref{SOURCE}) of class $L_p(D_*)$ in the domain
$D_*:=f(D)$ that satisfies the corresponding Neumann condition
(\ref{eqLIMHharm}).

\section{Poincar\'{e} problem in Bagemihl--Seidel system}

Arguing similarly to the last section, see also Section 4, we obtain
by Theorem 6 in \cite{GNRY} the following statement.

\medskip

{\bf Theorem 7.}{\it\, Let $D$ be a Jordan domain in $\mathbb C$,
$\nu:\partial D\to\mathbb C$, $|\nu (\zeta)|\equiv 1$, and
$\varphi:\partial D\to\mathbb C$ be measurable functions with
respect to the logarithmic capacity and let $\{
\gamma_{\zeta}\}_{\zeta\in\partial D}$ be a family of Jordan arcs of
class ${\cal{BS}}$ in ${D}$. Suppose that $A\in M^{2\times 2}(D)$
has H\"older continuous entries and $g:D\to\mathbb R$ is of class
$L_p(D)$ for some $p>2 $.

\medskip

Then there exist smooth weak solutions $u:D\to\mathbb R$ of equation
(\ref{6.8}) that have the limits along $\gamma_{\zeta}$
\begin{equation}\label{eqLIMITBS} \lim\limits_{z\to\zeta , z\in D}\ \frac{\partial u}{\partial \nu}\ (z)\ =\
\varphi(\zeta) \quad\quad\quad \mbox{q.e.\ on $\ \partial D$}\
.\end{equation} Furthermore, the space of such solutions $u$ has
infinite dimension.}

\bigskip

{\bf Remark 10.} Moreover, each such solution $u$ has the
representation $u =U\circ f|_D$, where $f=f^{\mu}:\mathbb
C\to\mathbb C$ is a quasiconformal mapping with $\mu$ in
(\ref{6.10}) extended onto $\mathbb C$ outside of $D$ as in the
proof of Theorem 6, and $U$ is a generalized harmonic function with
the source $G$ in (\ref{SOURCE}) of class $L_p(D_*)$ in the domain
$D_*:=f(D)$ that satisfies the corresponding Poincare condition
(\ref{eqLIMHharm}) along $\Gamma_{\xi}:=f(\gamma_{f^{-1}(\xi)}),
\xi\in\partial D_* $.


\section{Riemann--Poincar\'{e} type problems}

Arguing similarly to the last section and to Sections 5 and 6, we
reduce the proof of existence of solutions $u^{\pm}$ in the next
statement to Corollary 5 in \cite{GNRY} on such problems for the
corresponding generalized harmonic functions $U^{\pm}$ with the
source $G$ in (\ref{SOURCE}) through the representation
$u^{\pm}=U^{\pm}\circ f^{\mu}$, where
 $f=f^{\mu}:\mathbb C\to\mathbb C$ is a quasiconformal mapping with
$\mu$ in (\ref{6.10}).

\medskip

{\bf Theorem 8.}{\it\, Let $D$ be a Jordan domain in ${\mathbb C}$,
$\varphi :\partial D\times\mathbb R\to\mathbb R$ satisfy the
Carath\'{e}odory conditions, $\nu$ and $\nu_*:\partial D\to\mathbb
C$, $|\nu (\zeta)|\equiv 1$, $|\nu_* (\zeta)|\equiv 1$, are
measurable with respect to the logarithmic capacity, and let
$\{\gamma^+_{\zeta}\}_{\zeta\in\partial D}$ and
$\{\gamma^-_{\zeta}\}_{\zeta\in\partial D}$ be families of Jordan
arcs of class ${\cal{BS}}$ in ${D}$ and $\mathbb
C\setminus\overline{ D}$, correspondingly. Suppose that $A\in
M^{2\times 2}(\mathbb C)$ has H\"older continuous entries and
$g:\mathbb C\to\mathbb R$ is of class $L_p(\mathbb C)$ for some $p>2
$ with compact support.

\medskip

Then there exist smooth weak solutions $u^+:D\to\mathbb R$,
$u^-:\mathbb C\setminus\overline D\to\mathbb R$ of equation
(\ref{6.8}) such that
\begin{equation}\label{eqMIXED} \left[\frac{\partial u}{\partial\nu_*}\right]^+(\zeta)\ =\ \varphi\left(\,\zeta,\,
\left[\frac{\partial u}{\partial\nu}\right]^- (\zeta)\, \right)
\quad\quad\quad \mbox{q.e.\ on $\ \partial D$}\ ,
\end{equation}
where $\left[\frac{\partial u}{\partial\nu_*}\right]^+(\zeta)$ and
$\left[\frac{\partial u}{\partial\nu}\right]^-(\zeta)$ are limits of
the directional derivatives $\frac{\partial u^+}{\partial\,
\nu_*\,}\ (z)$ and $\frac{\partial u^-}{\partial\, \nu\,}\ (z)$ as
$z\to\zeta$ along $\gamma^+_{\zeta}$ and $ \gamma^-_{\zeta}$,
correspondingly.


Furthermore, the function $\left[\frac{\partial
u}{\partial\nu}\right]^- (\zeta)$ can be arbitrary measurable with
respect to the logarithmic capacity and, correspondingly, the space
of all such couples $(u^+,u^-)$ has the infinite dimension for any
such prescribed functions $A$, $\varphi$, $\nu$, $\nu_*$ and
collections $\gamma^+_{\zeta}$ and $ \gamma^-_{\zeta}$,
$\zeta\in\partial D$. }

\bigskip

Finally note that the present paper creates the basis for articles
on the corresponding results for semi-linear equations (with
nonlinear sources) of mathematical physics that will be published
elsewhere.

\bigskip

{\bf \noindent Vladimir Gutlyanskii, Olga Nesmelova, Vladimir Ryazanov} \\
Institute of Applied Mathematics and Mechanics\\ of National Academy
of Sciences of Ukraine,\\ Slavyansk 84 100,  UKRAINE\\
vgutlyanskii@gmail.com, star-o@ukr.net,\\ vl.ryazanov1@gmail.com,
Ryazanov@nas.gov.ua

\bigskip

\noindent {\bf Eduard Yakubov}\\
Holon Institute of Technology, Holon, Israel,\\
yakubov@hit.ac.il, eduardyakubov@gmail.com

\end{document}